\documentclass[reqno]{article}
\usepackage{graphicx}
\usepackage{color}
\usepackage{float}
\usepackage{amsmath}
\usepackage{cancel}

\newlength{\textlarg}
\newcommand{\barre}[1]{%
\settowidth{\textlarg}{#1}
#1\hspace{-\textlarg}%
\rule[0.5ex]{\textlarg}{1pt}}

\newcommand{\barred}[1]{%
\settowidth{\textlarg}{#1}
#1\hspace{-\textlarg }%
\hspace{-2pt}\rule[1ex]{\textlarg}{0.4pt}
\hspace{-\textlarg}\hspace{-3pt}%
\rule[0.1ex]{\textlarg}{0.4pt}}

\providecommand{\MSC}[1]{\textbf{\it{MSC:}} #1}
\providecommand{\keywords}[1]{\textbf{\it{Keywords:}} #1}
\providecommand{\Subject}[1]{\textbf{\it{Subject:}} #1}
\begin{document}

\pagestyle{myheadings}

\title{\bf Remarks on the Egyptian 2/D table in favor of a global approach
(D composite number) }
\author{\tt Lionel Br\'{e}hamet \\
Retired research scientist \\
\date{}
{\tt brehamet.l@orange.fr}}
\maketitle
\begin{flushleft}

\begin{abstract}

 \noindent
Egyptian decompositions of 2/D as a sum of two unit fractions are studied by means of certain divisors of D, namely r and s.
Our analysis does not concern the method to find r and s, but just why the scribes have chosen a solution instead of another.
A comprehensive approach, {\sf unconventional}, is developed which implies having an overview of  all pre-calculated alternatives.  Difference s-r is the basis of the classification to be examined before taking the right decisions, step by step, if difficulties are encountered.  An adequate adjustment of s-r  limits the number of all alternatives to 57, what is very few. 
A four-component generator (2/3, 2/5, 2/7, 2/11) operates as a (hidden) mother-table. Adding few logical rules of common sense is enough  to explain the reasons of the Egyptian choices. 
Even the singular case 2/95, which can not be decomposed into two fractions but only into three, reveals a doubly justified explanation. 
\\ \vspace{1em}

\noindent
\Subject{math.HO}\\
\noindent
\MSC{01A16}\\
\noindent

\end{abstract}

\newpage
\hspace*{1.5em}A general presentation of the Egyptian 2/D table and its bibliographic references have been done in our first paper on the subject
{\bf\cite {Brehamet1}}. A now famous scribe named Ahmes was the transcriber of the Rhind Mathematical Papyrus (RMP) which contains the table. \\
It should be kept in mind that the so-called '$2/D$ table' is not the work of only one scribe, but surely results of indeterminate periods of trials and improvements processed
by a team of scribes.\\ 
\hspace*{1.5em}In this second paper, we shall analyze only the decompositions into 2 terms if $D$ is composite. \\
Under a modern form we may write:
\begin{equation}
\frac{2}{D}= \frac{1}{\mathcal D_1}+ \frac{1}{\mathcal D_2}.
\label{eq:InitialProblem}
\end{equation}
Be noted that here the denominators of both unit fractions are written with calligraphic letters.

					\section{Presentation of the data  }
						\label{PART I}
\setcounter{section}{1}
\setcounter{equation}{0}

\hspace*{1.5em}Before we start choosing decompositions for composite numbers $D$, it can be assumed that in the ancient time, 
the results for the prime numbers were already obtained. This does not necessarily mean that there was a real chronology, like `results for $D$ prime' before
`2/$D$ table for $D$ composite'. These two aspects are interrelated and the project of the whole table was certainly been processed simultaneously, but in a  interdependent way by a few workshops of talented scribes. \\

Accordingly, as it was stressed in our preceding paper {\bf\cite{Brehamet1}}, we have at our disposal
a powerful tool, namely the four first simplest \sf{[two-terms]} \rm   decompositions of $2/D$, namely $2/3$, $2/5$, $2/7$ and $2/11$. 
Every composite up to $99$ (=$3\times 33$=$11\times 9$) may be derived from a  Mother-Table  looking like this:
\begin{table}[htp] 
\caption{\sf  Mother-Table (a tool for $D$ composite)}
\small
\begin{center}
\begin{tabular}{|l|}
 \hline
$2/3=1/2+1/6\,_ {\textcolor{red}{2}}$   \\ [0.01in] 
$2/5=1/3+1/15\,_{\textcolor{red}{ 3}}$   \\ [0.01in] 
$2/7=1/4+1/28\,_ {\textcolor{red}{4}}$   \\ [0.01in] 
$\cdots\cdots\cdots\cdots\cdots\cdots\cdot$  \\ [0.01in] 
$2/11=1/6+1/66\,_{\textcolor{red}{ 6}}$   \\ \hline
\end{tabular}
\end{center}
\label{MotherTable}
\end{table}\\

\hspace*{1.5em}This table is far from insignificant and reveals a lot of information.\\
\hspace*{1.5em}$\boxed{\sf\alpha\rm}$ Ability to have found all decompositions (unique)  of $2/D$ with $D$ prime.\\
\hspace*{1.5em}$\boxed{\sf\beta\rm}$ Deliberate choice of stopping [2-terms] decompositions from $11$. \\
\hspace*{1.5em}$\boxed{\sf\gamma\rm}$ As maximal multiplicity of the \it last denominator \rm with $D$ is {\sf 6}, keep this simple property \\
\hspace*{3.2em}for all the other remaining decompositions with $D$ composite
				\footnote{As will be seen further in all tables \ref{Complete2terms},  \ref{2termsconflict}, and  \ref{Nodd}, this decision is quite achievable.}.\\
\vspace{0.2em}
Below we reorder  the results  given by Ahmes, by using {\textcolor{red}{our $m_2$ numbers }}, multiplicity of ${\mathcal D_2}$   with $D$ (when there is one). Singular case $2/95$ [3-terms] is written separately and just aligned with $2/55$.
\\ \vspace{0.2em}
\begin{table}[htbp]
\caption{\sf  Ahmes's selections for D composite}
\small
\begin{center}
				\begin{tabular}{ll}
\begin{tabular}{|l|}\hline
$2/9=1/6+1/18\,_ {\textcolor{red}{2}}$   \\
$2/15=1/10+1/30\,_{\textcolor{red}{ 2}}$   \\ 
$2/21=1/14+1/42\,_{\textcolor{red}{ {2}}}$   \\ 
$2/25=1/15+1/75\,_{\textcolor{red}{ {3}}}$   \\ 
$2/27=1/18+1/54\,_{\textcolor{red}{ {2}}}$   \\  
$2/33=1/22+1/66\,_{\textcolor{red}{ {2}}}$   \\ 
$2/35=1/30+1/42  $  \\ 
$2/39=1/26+1/78\,_{\textcolor{red}{ {2}}}$   \\ 
$2/45=1/30+1/90\,_{\textcolor{red}{ {2}}}$   \\ 
$2/49=1/28+1/196\,_{\textcolor{red}{ {4}}}$   \\ 
$2/51=1/34+1/102\,_{\textcolor{red}{ {2}}}$   \\  \hline
\end{tabular}
&
\begin{tabular}{|l|}\hline
$2/57=1/38+1/114\,_{\textcolor{red}{ {2}}}$   \\
$2/63=1/42+1/126\,_{\textcolor{red}{ {2}}}$   \\
$2/65=1/39+1/195\,_{\textcolor{red}{ {3}}}$   \\ 
$2/69=1/46+1/138\,_{\textcolor{red}{ {2}}}$   \\ 
$2/75=1/50+1/150\,_{\textcolor{red}{ {2}}}$   \\
$2/77=1/44+1/308\,_{\textcolor{red}{ {4}}}$   \\  
$2/81=1/54+1/162\,_{\textcolor{red}{ {2}}}$   \\ 
$2/85=1/51+1/255\,_{\textcolor{red}{ {3}}}$   \\ 
$2/87=1/58+1/174\,_{\textcolor{red}{ {2}}}$   \\ 
$2/91=1/70+1/130  $  \\  
$2/93=1/62+1/186\,_{\textcolor{red}{ {2}}}$   \\ 
$2/99=1/66+1/198\,_{\textcolor{red}{ {2}}}$   \\ \hline
\end{tabular}
				\end{tabular}
\end{center}
\begin{center}
	\begin{tabular}{ll}
\begin{tabular}{|l|}\hline
$2/55=1/30+1/330\,_{\textcolor{red}{ {6}}}$   \\  \hline
\end{tabular}
					&
\begin{tabular}{|l|}\hline
$2/95=1/60+1/380\,_{\textcolor{red}{ 4}}+1/570\,_{\textcolor{red}{6}}$   \\  \hline
\end{tabular}
	\end{tabular}				
\end{center}

\label{Composites}
\end{table}

					\section{All possible solutions (a finite countable set)}
						\label{PART II}
\setcounter{section}{2}
\setcounter{equation}{0}

\renewcommand{\theequation }{\Roman{section}.\arabic{equation}}

\hspace{1.5em}As soon the matter of composite number is approached, it becomes ineluctable to consider how the number can be decomposed in products of its divisors,
as well for ancient mathematicians as for modern like us. 
It is possible  to calculate all solutions for $D$ composite, by  using a theorem little-known or ignored. We call it {\mathversion{bold}$N^{[2]}\, \top \mbox{\tt heorem }$},
whose we will detail some aspects.
It is unlikely that these properties were known to the ancients, but an heuristic discovery is never to exclude.\\ 
\vspace{0.5em}
\hspace*{7em}Equation to be solved:
\begin{equation}
N^2-\sigma_1 N+\sigma_2 =0 , \mbox{\hspace{1em}with integral solutions $N_1$, $N_2 >N_1$}.
\end{equation}
\hspace*{7em}Condition to be filled:
\begin{equation}
\frac{\sigma_1}{\sigma_2 }=  \frac {p}{q}\mbox{\hspace{1em}= \hspace{0.2em}irreducible fraction $< 1$ }.
\end{equation}
\hspace*{7em}Method: decompose $q$ in triplets of its divisors as follows
\begin{equation}
q= k \times r \times s \mbox{\hspace{1em} \it with }r<s .
\label{eq:Decqinkrs}
\end{equation}
\hspace*{7em}All the solutions  are found in the finite set of $N_1$, $N_2$ given by 
\begin{equation}
N_1= \frac{k r (r+s)}{p} \mbox{\hspace{1em} and \hspace{1em}}N_2= \frac{k s(r+s)}{p},
\label{eq:SolutionsN1N2}
\end{equation}
	
\begin{equation}
\mbox{where $(r+s)$ is selected such that }
\left(\frac{r+s}{p}\right) \;\mbox{to be integral.}
						\footnote{Of course selection due to  Eq.(\ref{eq:pSelection}) vanishes for the Egyptian table since $p=2$, and $q$, $r$, $s$ are odd. } 
\label{eq:pSelection}
\end{equation}
\begin{equation}
\mbox{\sf Egyptian property:}\quad \frac{1}{N_1} + \frac{1}{N_2}= \frac{p}{q}.
\end{equation}
The relation between both solutions $N_1,N_2$ depends only on $r,s$:
\begin{equation}
 \frac{N_2}{N_1}=\frac{s}{r}.
\end{equation}
\mathversion{normal}
The ratio between $N_2$ and $q$ is given by
\begin{equation}
\frac{N_2}{q}=\frac{r+s}{pr}.
\end{equation}
If this quantity is integral, it will  be labeled  conventionally by {\sf $m_2$} (a true multiplicity), 
if not we will use the denomination `fractional multiplicity' with the label $\mu_2$.\\
Another  important definition will be useful, namely
\begin{equation}
\Delta_{\mathbf {r}}^{\mathbf {s}}\mathbf{=(s-r)}.
\end{equation}
It is the gap between two divisors of $D$ which measures the closeness between them.\\
Of course, the first trivial triplet decomposition of $q=1 \times 1 \times q$ represents nothing but a trivial identity, 
since the integral solutions $\displaystyle N_1=\frac{q+1}{p}$, $\displaystyle N_2 =\frac{q(q+1)}{p}$ satisfy
\begin{equation}
\frac{1}{N_1}+ \frac{1}{N_2}= \frac{p}{q+1}+ \frac{p}{q(q+1)}= \frac{p}{q}.
\label{eq:Primeformula}
\end{equation}
In the Egyptian $2/q$ table ($p=2$), if $q$ is prime, formula (\ref{eq:Primeformula}) yields the unique solution. Clearly found by the scribes.
However for $q$ composite, this remains a solution among other.\\
This extremal  solution produces the largest denominator $N_2$ and the highest multiplicity of $N_2$ with $q$, namely $m_2=(q+1)/2$, 
certainly discarded by Egyptian scribes for obvious reasons. \\ 
Furthermore there are other solutions to be discarded also according to our comment $\boxed{\sf\gamma\rm}$. \\
This is consistent with the last value  of the Mother-table, namely $11$.
Indeed, if we admit a limitation of $m_2$ by a top-flag like
\begin{equation}
\mbox{\boldmath $\top$}\!\!_ f^{\;\;[2]}= {\mathbf 6},
\end{equation}
then, for $p=2$, the following inequality is found:
\begin{equation}
s\leq \mathbf{11 }\,r. 
\label{eq:Number11}
\end{equation}
Yes, number $\mathbf{11}$ is well the one that appears in our comment  $\boxed{\sf\beta\rm}$.\\
It has often been debated about this value $11$ and the Eratosthenes's sieve as well as about the ``perfectness'' of the number $6$, see for example Brown {\bf \cite{Brown}},
however these notions are quite useless here. Simply, the present approach brings up an unexpected link between the numbers 11 and 6.\\
Of course inequality (\ref{eq:Number11}) seriously reduces the number of trials for finding convenient integers $r$, $s$.
In our table the solutions excluded by this inequality will be displayed as lines doubly strikethrough, like $2/D=$\barred{$1/\mathcal D_1+1/\mathcal D_2$}.\\ 
\hspace*{1.5em}Be aware that the authors of the table have succeeded in calculating $2/35$ or $2/91$ where there is no  multiple  of $D$, 
then the question is not `why these two singular cases have been chosen?' but why other existing cases have been discarded?\\
\hspace*{1.5em}Since they have also certainly succeeded to calculate all possible options, we are going to list all the cases ordered 
according to increasing values of $\Delta_{\mathbf {r}}^{\mathbf {s}}$. The method is quite similar to  that followed 
in our preceding analysis {\bf \cite{Brehamet1}}, where it can be observed that all Egyptian  solutions are found in the set \{$\Delta_d, \Delta_{d'}\leq 10$\}.
{\sf The value 10 }(Egyptian decimal stem)  obviously is not a universal constant, but seems {\sf unavoidable.} Other general problems of arithmetic 
would probably not have benefited from such an opportunity. However once again it can be tried to limit the investigations by the following arithmetical prescription:
\begin{equation}
\Delta_{\mathbf {r}}^{\mathbf {s}}\leq 10.
\end{equation}
We shall use a slash to indicate that an element  does not satisfy this constraint, like $\cancel{\Delta_{r}^{s}}$.\\
Condition $\Delta_{\mathbf {r}}^{\mathbf {s}}\leq 10$ is an `absolute' condition stronger than the relative issued 
from (\ref{eq:Number11}).\\
\hspace*{1.5em}{\sf Although here the divisors involved ($s$ and $r$) are not those of the first denominator  as in {\bf \cite{Brehamet1}} but  those of $D$, we were surprised 
that the same approach was still operating and fruitful.}\\ Moreover, as we shall see, the  "intriguing" cases {\sf 2/35} and {\sf 2/91} do not need to have a separate treatment, as invoke a ``arithmetic-harmonic'' decomposition  {\bf \cite{Brown}}, (post) dating from ancient Greeks.\\ 
\hspace*{1.5em}Considering such examples (and some other) as singular, or ``{\sl not optimal }''  {\bf \cite {GardnerMilo}}, reflects  just  a modern point of view, mathematically speaking.
An optimality concept depends on  the views, and was certainly  different in ancient context.\\
{\sf Note:} for interested readers, we have added two appendices ({\tt A} and {\tt B}) where we
sketch out two different methods leading to the necessity of decomposing $D$ or $D^2$ as a product of two of its divisors. \\
Total cases to be analyzed (as reported in our abstract) does not take into account the cases involving a slash or lines doubly strikethrough. 
Throughout this paper, Egyptian choices will be denoted by the symbol \scriptsize ${Eg}$\normalsize. An overview of all possibilities is displayed
in Table \ref{Complete2terms}.

\begin{table}[htbp]
\caption{\sf 2 terms decompositions with increasing values of $ \Delta_{r}^{s}$   }
\begin{center}
\scriptsize
							\begin{tabular}{ll}
	
\begin{tabular}{|c|l|}\hline
$ \Delta_{r}^{s}$   & \sf All possible solutions \\ \hline \hline
$\boxed{2}$&$2/9=1/6+1/18\,_ {\textcolor{red}{2}}$\scriptsize $\;\,{Eg}\normalsize$  \\ 
\hline \hline
$\boxed{2}$&$2/15=1/10+1/30\,_{\textcolor{red}{ 2}}$\scriptsize $\;\,{Eg}\normalsize$  \\ 
2&$2/15=1/12+1/20\,_{\textcolor{red}{4/3}}$   \\ \hline \hline
$\boxed{2}$&$2/21=1/14+1/42\,_{\textcolor{red}{ {2}}}$\scriptsize $\;\,{Eg}\normalsize$   \\ \hline \hline
$\boxed{2}$&$2/27=1/18+1/54\,_{\textcolor{red}{ {2}}}$\scriptsize $\;\,{Eg}\normalsize$  \\  \hline  \hline
$\boxed{2}$&$2/33=1/22+1/66\,_{\textcolor{red}{ {2}}}$\scriptsize $\;\,{Eg}\normalsize$   \\ 
\hline \hline
$\boxed{2}$&$2/35=1/30+1/42\,_{\textcolor{red}{6/5 }}  $\scriptsize $\;\,{Eg}\normalsize$  \\   
\hline \hline 
$\boxed{2}$&$2/39=1/26+1/78\,_{\textcolor{red}{ {2}}}$\scriptsize $\;\,{Eg}\normalsize$   \\ 
\hline \hline
$\boxed{2}$ &$2/45=1/30+1/90\,_{\textcolor{red}{ {2}}}$\scriptsize $\;\,{Eg}\normalsize$   \\ 
2&$2/45=1/36+1/60\,_{\textcolor{red}{4/3 }}  $  \\ \hline \hline 
$\boxed{2}$&$2/51=1/34+1/102\,_{\textcolor{red}{ {2}}}$\scriptsize $\;\,{Eg}\normalsize$   \\ 
\hline \hline
$\boxed{2}$&$2/57=1/38+1/114\,_{\textcolor{red}{ {2}}}$\scriptsize $\;\,{Eg}\normalsize$   \\ \hline \hline
$\boxed{2}$&$2/63=1/42+1/126\,_{\textcolor{red}{ {2}}}$\scriptsize $\;\,{Eg}\normalsize$   \\ 
2&$2/63=1/56+1/72\,_{\textcolor{red}{8/7 }}  $  \\  \hline \hline
$\boxed{2}$&$2/69=1/46+1/138\,_{\textcolor{red}{ {2}}}$\scriptsize $\;\,{Eg}\normalsize$   \\ 
\hline \hline
$\boxed{2}$&$2/75\mathit{_a}=1/50+1/150\,_{\textcolor{red}{ {2}}}$\scriptsize $\;\,{Eg}\normalsize$   \\ 
2&$2/75=1/60+1/100\,_{\textcolor{red}{4/3}}  $ \\ 
\hline \hline
$\boxed{2}$ &$2/81\mathit{_a}=1/54+1/162\,_{\textcolor{red}{ {2}}}$\scriptsize $\;\,{Eg}\normalsize$   \\ 
\hline 
$\boxed{2}$&$2/87=1/58+1/174\,_{\textcolor{red}{ {2}}}$\scriptsize $\;\,{Eg}\normalsize$   \\ 
\hline \hline
$\boxed{2}$&$2/93=1/62+1/186\,_{\textcolor{red}{ {2}}}$\scriptsize $\;\,{Eg}\normalsize$   \\ 
\hline \hline
$\boxed{2}$&$2/99=1/66+1/198\,_{\textcolor{red}{ {2}}}$\scriptsize $\;\,{Eg}\normalsize$   \\ 
2&$2/99=1/90+1/110\,_{\textcolor{red}{10/9 }}  $  \\ \hline 
\hline \hline
4&$2/15=1/9+1/45\,_{\textcolor{red}{ 3}}$  \\  \hline \hline 
4&$2/21=1/15+1/35\,_{\textcolor{red}{5/3}}$  \\ 
\hline \hline
$\boxed{4}$&$2/25=1/15+1/75\,_{\textcolor{red}{ {3}}}$\scriptsize $\;\,{Eg}\normalsize$ \\ 
\hline \hline
4&$2/35=1/21+1/105\,_{\textcolor{red}{ {3}}}$\\   \hline \hline
4&$2/45\mathit{_a}=1/27+1/135\,_{\textcolor{red}{ {3}}}$  \\ \hline \hline
4&$2/55=1/33+1/165\,_{\textcolor{red}{ {3}}}$\\   \hline \hline
4&$2/63=1/45+1/105\,_{\textcolor{red}{5/3}}  $ \\  \hline  \hline
$\boxed{4}$&$2/65=1/39+1/195\,_{\textcolor{red}{3 }}  $\scriptsize $\;\,{Eg}\normalsize$  \\   \hline \hline
4&$2/77=1/63+1/99\,_{\textcolor{red}{9/7 }}  $\\   
\hline \hline
$\boxed{4}$&$2/85=1/51+1/255\,_{\textcolor{red}{3}}  $\scriptsize $\;\,{Eg}\normalsize$  \\  \hline \hline 
4&$2/95\stackrel{\mathbf?}{=}1/57+1/285\,_{\textcolor{red}{ {3}}}$ \\   \hline \hline \hline
6&$2/21=1/12+1/84\,_{\textcolor{red}{ {4}}}$   \\ \hline \hline
6&$2/27\mathit{_a}=1/15+1/135\,_{\textcolor{red}{ {5}}}$  \\ \hline\hline
6&$2/35=1/20+1/140\,_{\textcolor{red}{ {4}}}$   \\   \hline \hline
$\boxed{6}$ &$2/49=1/28+1/196\,_{\textcolor{red}{ {4}}}$\scriptsize $\;\,{Eg}\normalsize$   \\ \hline \hline
6&$2/55=1/40+1/88\,_{\textcolor{red}{8/5}}$   \\   \hline \hline
6&$2/63\mathit{_a}=1/36+1/252\,_{\textcolor{red}{ {4}}}$   \\ \hline \hline 
$\boxed{6}$& $2/77=1/44+1/308\,_{\textcolor{red}{4 }} $\scriptsize $\;\,{Eg}\normalsize $ \\ \hline \hline
6&$2/81\mathit{_b}=1/54+1/162\,_{\textcolor{red}{ {2}}}$  \\ \hline \hline 
6 &$2/91=1/52+1/364\,_{\textcolor{red}{ {4}}}$ \\  
$\boxed{6}$&$2/91=1/70+1/130\,_{\textcolor{red} {10/7}}$\scriptsize $\;\,{Eg}\normalsize$   \\   
 \hline \hline \hline
8&$2/27\mathit{_b}=1/15+1/135\,_{\textcolor{red}{ {5}}}$  \\ \hline \hline 
8&$2/33=1/21+1/77\,_{\textcolor{red}{7/3 }}  $ \\ \hline \hline 
8&$2/45=1/25+1/225\,_{\textcolor{red}{ {5}}}$  \\ \hline \hline 
8&$2/63=1/35+1/315\,_{\textcolor{red}{ {5}}}$   \\  \hline \hline 
8&$2/65=1/45+1/117\,_{\textcolor{red}{ 9/5}  }$  \\  \hline \hline 
8&$2/81\mathit{_{a'}}=1/45+1/405\,_{\textcolor{red}{ {5}}}$  \\ \hline \hline 
8&$2/99=1/55+1/495\,_{\textcolor{red}{ {5}}}$  \\ \hline \hline 
8&$2/99=1/63+1/231\,_{\textcolor{red}{7/3}}$  \\ \hline \hline \hline 
10&$2/33=1/18+1/198\,_{\textcolor{red}{ {6}}}$   \\ \hline \hline 
10&$2/39=1/24+1/104\,_{\textcolor{red}{8/3 }}$  \\ \hline  
\end{tabular}
								&
\begin{tabular}{|c|l|}\hline
$ \Delta_{r}^{s}$   & \sf All possible solutions \\ \hline \hline
$\boxed{\mathbf {10}}$&$\mathbf{2/55=1/30+1/330\,_{\textcolor{red}{ {6}}}}$\scriptsize $\;\,{\bf Eg}\normalsize$   \\   \hline \hline
10&$2/75\mathit{_b}=1/50+1/150\,_{\textcolor{red}{ {2}}}$\\ \hline \hline 
10&$2/77=1/42+1/462\,_{\textcolor{red}{ {6}}}$   \\   \hline \hline
10&$2/99\mathit{_a}=1/54+1/594\,_{\textcolor{red}{ {6}}}$   \\ \hline \hline \hline
\cancel{12}&$2/39=$\barred{$1/21+1/273\,_{\textcolor{red}{ {7}}}$}  \\ \hline \hline 
\cancel{12}&$2/45\mathit{_b}=1/27+1/135\,_{\textcolor{red}{ {3}}}$  \\ \hline \hline 
\cancel{12}&$2/65=$\barred{$1/35+1/455\,_{\textcolor{red}{ {7}}}$}   \\   \hline \hline
\cancel{12}&$2/85=1/55+1/187\,_{\textcolor{red}{11/5 }} $ \\ \hline \hline 
\cancel{12}&  $2/91=$\barred{$1/49+1/637\,_{\textcolor{red}{ {7}}}$  }\\   \hline \hline \hline
\cancel{14}&$2/15=$\barred{$\!1/8+1/120\,_{\textcolor{red}{ 8}}$}  \\ \hline \hline 
\cancel{14}&$2/45=$\barred{$1/24+1/360\,_{\textcolor{red}{ {8}}}$}    \\ \hline \hline 
\cancel{14}&$2/51=1/30+1/170\,_{\textcolor{red}{10/3} }$  \\ \hline \hline 
\cancel{14}&$2/75=$\barred{$1/40+1/600\,_{\textcolor{red}{ {8}}}$}   \\ \hline \hline  
\cancel{14} &$2/95\stackrel{\mathbf?}{=}1/60+1/228\,_{\textcolor{red}{ 12/5}}$   \\ \hline \hline
\cancel{16}&$2/51=$\barred{$1/27+1/459\,_{\textcolor{red}{ {9}}}$}   \\ \hline \hline 
\cancel{16}&$2/57=1/33+1/209\,_{\textcolor{red}{11/3 }}$  \\ \hline \hline 
\cancel{16}&$2/85=$\barred{$1/45+1/765\,_{\textcolor{red}{ {9}}}$}   \\ \hline \hline \hline
\cancel{18}&$2/57=$\barred{$1/30+1/570\,_{\textcolor{red}{ {10}}}$}   \\ \hline \hline 
\cancel{18}&$2/63\mathit{_b}=1/36+1/252\,_{\textcolor{red}{ {4}}}$   \\ \hline \hline \hline
\cancel{20}&$2/21=$\barred{$1/11+1/231\,_{\textcolor{red}{ {11}}}$}  \\ \hline \hline 
\cancel{20}&$2/63=$\barred{$1/33+1/693\,_{\textcolor{red}{ {11}}}$}   \\  \hline \hline \hline 
\cancel{22}&$2/69=$\barred{$1/36+1/828\,_{\textcolor{red}{ {12}}}$}   \\  \hline \hline 
\cancel{22}&$2/75=1/42+1/350\,_{\textcolor{red}{14/3 }}$  \\ \hline \hline \hline 
\cancel{24}&$2/25=$\barred{$1/13+1/325\,_{\textcolor{red}{ {13}}}$}\\ \hline \hline 
\cancel{24}&$2/75=$\barred{$1/39+1/975\,_{\textcolor{red}{ {13}}}$}   \\ \hline \hline 
\cancel{24}&$2/81\mathit{_{b'}}=1/45+1/405\,_{\textcolor{red}{ {5}}}$  \\ \hline \hline \hline 
\cancel{26}&$2/27=$\barred{$1/14+1/378\,_{\textcolor{red}{ {14}}}$}  \\ \hline \hline 
\cancel{26}&$2/81=$\barred{$1/42+1/1134\,_{\textcolor{red}{ {14}}}$}   \\ \hline \hline 
\cancel{26}&$2/87=1/48+1/464\,_{\textcolor{red}{16/3}}$  \\ \hline \hline \hline 
\cancel{28}&$2/87=$\barred{$1/45+1/1305\,_{\textcolor{red}{ {15}}}$}   \\ \hline \hline 
\cancel{28}&$2/93=1/51+1/527\,_{\textcolor{red}{17/3}} $  \\ \hline \hline \hline 
\cancel{30}&$2/93=$\barred{$1/48+1/1488\,_{\textcolor{red}{ {16}}}$}  \\ \hline \hline 
\cancel{30}&$2/99\mathit{_b}=1/54+1/594\,_{\textcolor{red}{ {6}}}$   \\ \hline \hline \hline 
\cancel{32}&$2/33=$\barred{$1/17+1/561\,_{\textcolor{red}{ {17}}}$}   \\ \hline \hline  
\cancel{32}&$2/99=$\barred{$1/51+1/1683\,_{\textcolor{red}{ {17}}}$}   \\ \hline \hline \hline 
\cancel{34}&$2/35=$\barred{$1/18+1/630\,_{\textcolor{red}{ {18}}}$}   \\  \hline \hline \hline
\cancel{38}&$2/39=$\barred{$1/20+1/780\,_{\textcolor{red}{ {20}}}$}   \\ \hline \hline \hline 
\cancel{44}&$2/45=$\barred{$1/23+1/1035\,_{\textcolor{red}{ {23}}}$}   \\ \hline \hline \hline
\cancel{48}&$2/49=$\barred{$1/25+1/1225\,_{\textcolor{red}{ {25}}}$}  \\ \hline \hline \hline 
\cancel{50}&$2/51=$\barred{$1/26+1/1326\,_{\textcolor{red}{ {26}}}$}   \\ \hline \hline \hline 
\cancel{54}&$2/55=$\barred{$1/28+1/1540\,_{\textcolor{red}{ {28}}}$}   \\   \hline \hline \hline
\cancel{56}&$2/57=$\barred{$1/29+1/1653\,_{\textcolor{red}{ {29}}}$}   \\ \hline \hline \hline 
\cancel{62}&$2/63=$\barred{$1/32+1/2016\,_{\textcolor{red}{ {32}}}$}    \\ \hline \hline \hline 
\cancel{64}&$2/65=$\barred{$1/33+1/2145\,_{\textcolor{red}{ {33}}}$}  \\  \hline \hline \hline
\cancel{68}&$2/69=$\barred{$1/35+1/2415\,_{\textcolor{red}{ {35}}}$}   \\  \hline \hline \hline 
\cancel{74}&$2/75=$\barred{$1/38+1/2850\,_{\textcolor{red}{ {38}}}$}   \\ \hline \hline \hline 
\cancel{76}&$2/77=$\barred{$1/39+1/3003\,_{\textcolor{red}{ {39}}}$}   \\   \hline \hline  \hline
\cancel{80}&$2/81=$\barred{$1/41+1/3321\,_{\textcolor{red}{ {41}}}$}  \\ \hline \hline \hline 
\cancel{84}&$2/85=$\barred{$1/43+1/3655\,_{\textcolor{red}{ {43}}}$}   \\ \hline \hline \hline
\cancel{86}&$2/87=$\barred{$1/44+1/3828\,_{\textcolor{red}{ {44}}}$} \\ \hline  \hline \hline 
\cancel{90}&  $2/91=$\barred{$1/46+1/4186\,_{\textcolor{red}{ {46}}}$  }\\  \hline \hline \hline 
\cancel{92}&$2/93=$\barred{$1/47+1/4371\,_{\textcolor{red}{ {47}}}$}   \\ \hline \hline \hline 
\cancel{94}&$2/95=$\barred{$1/50+1/950\,_{\textcolor{red}{ {10}}}$} \\   \hline \hline \hline
\cancel{98}&$2/99=$\barred{$1/50+1/4950\,_{\textcolor{red}{ {50}}}$}   \\ \hline
\end{tabular}
							\end{tabular}
\end{center}

\label{Complete2terms}
\end{table}
\newpage
\hspace*{1.5em} We see that {\sf the  majority of the Egyptian choices were done according to the smallest values of $ \Delta_{r}^{s}$.} 
Nevertheless some cases show a `conflict' when two identical $ \Delta_{r}^{s}$  are in the presence of one-another.
As we will see below, this is not a hard problem to solve. 
Abdulaziz {\bf \cite {Abdulaziz}} is remained faithful to the spirit of ancient Egypt in substantiating  its arguments with fractional quantities, then we are going to do the same by using here our  `fractional multiplicity' $\mu_2$, however to a lesser extent. \\
The closeness between  $2$ and the `fractional multiplicity' $\mu_2$ will be estimated by a difference:\\
\begin{equation}
\mathit \Delta_\mu = 2 - \mu_2 ,\quad (\mbox{this $\mathit \Delta$ here is written in \it italic}).
\end{equation}
This difference is also fractional and could not have been a source of difficulty to sort them in increasing order, for example. In this kind of work, the ancient scribes were talented!\\
\begin{table}[htbp]
\caption{\sf Conflicting cases with identical $ \Delta_{r}^{s}$ }
\begin{center} 
\begin{tabular}{|c|l|l|l|}\hline
$ \Delta_{r}^{s}$   & \sf All possible solutions & Appreciations & \sf Decisions \\ \hline \hline
$\boxed{2}$&$2/15=1/10+1/30\,_{\textcolor{red}{ 2}}$\scriptsize $\;\,{Eg}\normalsize$  & $m_2=2$ possible & $m_2=2$\\ 
2&$2/15=1/12+1/20\,_{\textcolor{red}{\frac{4}{3} }}$  & $ \mathit \Delta_\mu = \frac{42}{63}$ &\\ \hline \hline 
$\boxed{2}$ &$2/45=1/30+1/90\,_{\textcolor{red}{ {2}}}$\scriptsize $\;\,{Eg}\normalsize$  &$m_2=2$ possible & $m_2=2$\\ 
2&$2/45=1/36+1/60\,_{\textcolor{red}{\frac{4}{3} }}  $  &$ \mathit \Delta_\mu =\frac{42}{63}$ &\\ \hline \hline  
$\boxed{2}$&$2/63=1/42+1/126\,_{\textcolor{red}{ {2}}}$\scriptsize $\;\,{Eg}\normalsize$  &$m_2=2$ possible & $m_2=2$\\ 
2&$2/63=1/56+1/72\,_{\textcolor{red}{\frac{8}{7} }}  $ & $ \mathit \Delta_\mu =\frac{54}{63}$ &\\  \hline \hline 
$\boxed{2}$&$2/75\mathit{_a}=1/50+1/150\,_{\textcolor{red}{ {2}}}$\scriptsize $\;\,{Eg}\normalsize$ & $m_2=2$ possible & $m_2=2$\\ 
2&$2/75=1/60+1/100\,_{\textcolor{red}{\frac{4}{3} }}  $  &$ \mathit \Delta_\mu =\frac{42}{63}$ & \\ \hline \hline 
$\boxed{2}$&$2/99=1/66+1/198\,_{\textcolor{red}{ {2}}}$\scriptsize $\;\,{Eg}\normalsize$   &$m_2=2$ possible & $m_2=2$\\ 
2&$2/99=1/90+1/110\,_{\textcolor{red}{\frac{10}{9} }}  $ &$ \mathit \Delta_\mu =\frac{56}{63}$ &\\ \hline \hline \hline
6 &$2/91=1/52+1/364\,_{\textcolor{red}{ {4}}}$ & $m_2=4$ too far from $2$ &\\ 
$\!\!\cancel{\exists}$\;6 :&$2/91=1/{\mathcal D_1}+1/{\mathcal D_{2}}\,_{\textcolor{red}{ {2}}}$&virtual row: $\cancel{\exists}\;m_2 =2\;\,^{\ddagger}$  &\\ 
$\boxed{6}$&$2/91=1/70+1/130\,_{\textcolor{red}{ {\frac{10}{7}}}}$\scriptsize $\;\,{Eg}\normalsize$  &$ \mathit \Delta_\mu =\frac{36}{63}$
 & $\mu_2=\frac{10}{7}$\\   \hline 
\end{tabular}\\ 
\end{center}
\label{2termsconflict}
\end{table}
$^{\ddagger}$ \scriptsize in contrast to case 
\begin{tabular}{|c|} \hline 
$2/81\mathit{_b}=1/54+1/162\,_{\textcolor{red}{ {2}}}$ \\  \hline\end{tabular}
where also $ \Delta_{r}^{s}=6$. Our additional second row for 2/91 is obviously virtual, but useful for the clarity of our presentation. \normalsize\\ \vspace{0.5em}
Remark that $2/91$ is the only element that contains twice the scribal annotation  ``\sl find''  \rm on the Rhind Papyrus {\bf \cite {Abdulaziz}}. 
Particular attention was really brought to this matter. It is not a coincidence that it appears in Table \ref{2termsconflict}.
In our view, it is not comparable with $2/35$ which was solved immediately through its $ \Delta_{r}^{s}=2$, by following our classification which discards  $ \Delta_{r}^{s}=4,6$. 
Accordingly, as in our first paper {\bf \cite{Brehamet1}}, we stay faithful to the logic of Occam's razor: simplicity and minimal hypotheses. \\{\sl Ie}:
{\sf if $m_2=2$ is available, then retain it, 
or else adopt   a value of $\mu_2$  closest to $2$.}\\ \vspace{0.5em}
\hspace*{1.5em}Rarely a major project unfolds linearly, ie without a hitch. Nothing out of the ordinary, even today. Indeed there are still three cases unanswered, 
namely 2/77, 2/95 (and 2/55).\\
\begin{table}[htbp]
\caption{\sl`no odd denominators' \sf precept}
\begin{center}
								\begin{tabular}{ll}
\begin{tabular}{|c|l|}\hline
$ \Delta_{r}^{s}$   & \sf All possible solutions \\ \hline \hline

4&$\barre{2/77}=\barre{\;\,1/63+1/99}\;_{\textcolor{red}{\barre{\!9/7 }}}  $\\   \hline \hline 
$\boxed{6}$& $2/77=1/44+1/308\,_{\textcolor{red}{4 }} $\scriptsize $\;\,{Eg}\normalsize $ \\ \hline \hline
10&$2/77=1/42+1/462\,_{\textcolor{red}{ {6}}}$   \\   \hline 

\end{tabular}
									&
\begin{tabular}{|c|l|}\hline
$ \Delta_{r}^{s}$   & \sf All possible solutions \\ \hline \hline
4&$\barre{2/95}\stackrel{\mathbf?}{=}\barre{\;\,1/57+1/285}\,_{\textcolor{red}{ {\barre{3}}}}$ \\   \hline

\end{tabular}
								\end{tabular}
\end{center}

\label{Nodd}
\end{table}
The `conflicting' cases have been solved and `forgotten'. Now another strategy must be applied in a independent way.
The aim is to solve one after the other difficulties encountered during the dynamical advancement of project. Thus a `local' decision should not interfere 
with those taken previously. In the present instance, a reasonable option appears to be that of the famous rule\\ {\sl`no odd denominators'}.
Only once we were forced to apply this option for solving the case  {\sf 2/89} into 4 terms {\bf \cite{Brehamet1}}. This was often invoked by Gillings {\bf\cite{Gillings}}, criticized by Bruins {\bf\cite{Bruins}} and nevertheless several times  used by Abdulaziz {\bf\cite{Abdulaziz}}. 
In accordance with this `precept' ({\sl here used two times}), Table \ref{Nodd} is displayed with some barred rows indicating what possibilities will be deleted. 
As will be shown further, there is no need apply this precept to {\sf 2/55}. There are two immediate effects:  {\sf the  smallest $ \Delta_{r}^{s}=6$ is chosen 
for  {\sf 2/77}} and there is no more possibility for  {\sf 2/95}! Thus we can write:\\ \vspace{0.5em}
\hspace*{1.5em}
 \begin{tabular}{|c|} \hline 
$2/77=1/44+1/308\,_{\textcolor{red}{4 }} $\scriptsize $\;\,{Eg}\normalsize $ \\  \hline\end{tabular} $\quad$ and $\quad$
\begin{tabular}{|c|} \hline 
$2/95 \; \mbox{\sf has no decomposition  into two terms}$ \\  \hline\end{tabular}\\ \vspace{0.5em}
Up to this point, the utilization rates of mother-table  \ref {MotherTable} with  various divisors,  are the following:\\ \vspace{0.5em}
\begin{tabular}{lll}
\begin{tabular}{|c|}
\hline
\sf Mother-rows \\ [0.01in] \hline
$2/3=1/2+1/6\,_ {\textcolor{red}{2}}$   \\ [0.01in] \hline
$2/5=1/3+1/15\,_{\textcolor{red}{ 3}}$   \\ [0.01in]  \hline
$2/7=1/4+1/28\,_ {\textcolor{red}{4}}$   \\ [0.01in]  \hline
$2/11=1/6+1/66\,_{\textcolor{red}{ 6}}$   \\ \hline
\end{tabular}
&
\begin{tabular}{|c|}
\hline
\it Number of times used \\ [0.01in] \hline
$16$   \\ [0.01in] \hline
$3$   \\ [0.01in]  \hline
$2$   \\ [0.01in]  \hline
$0$   \\ [0.01in]  \hline
\end{tabular}
& 
\begin{tabular}{|c|c|}\hline
\sf Examples  &  \sf Divisors\\ [0.01in] \hline
2/51 = [row$_{1}]/17$ & 17 \\ [0.01in] \hline
2/65 = [row$_{2}]/13$ & 13 \\ [0.01in] \hline
2/77 = [row$_{3}]/11$ & 11 \\ [0.01in] \hline
no =  [row$_{4}]/?? $ & ??
 \\ \hline
\end{tabular}
							\end{tabular}\\ \vspace{0.5em}

Now what becomes the case $2/55$? \\
\begin{table}[htbp]
\begin{center}
\begin{tabular}{|c|l|}\hline
$ \Delta_{r}^{s}$   & \sf Remaining solutions for 2/55\\ \hline \hline
4&$2/55=1/33+1/165\,_{\textcolor{red}{3}}$\\   \hline \hline
6&$2/55=1/40+1/88\,_{\textcolor{red}{8/5}}$   \\   \hline \hline
$\boxed{\mathbf {10}}$&$\mathbf{2/55=1/30+1/330\,_{\textcolor{red}{ {6}}}}$\scriptsize $\;\,{\bf Eg}\normalsize$   \\   \hline 
\end{tabular}
\end{center}
\end{table}
Please note that the relation of  $\sf 2/55$ with the mother-table is quite unique.
In our previous paper {\bf \cite {Brehamet1}} we showed that there were reasons for leaving the decompositions into  2 terms from the value $11$. This explained why start the remaining cases
 for $D$ prime in 
[3- or 4-terms] from $13$ (except for $23$ that we have quite justified). 
Furthermore this fact  assigned the quality of mother-table to the first four fractions. 
\sf If the case 2/55  does not  use
\{$2/11$\} as template, this quality is lost, and a lot of things become un-understandable or inconsistent. \rm This is why scribes adopted an \sf ultimate decision for 2/55\rm: \\

\begin{table}[htp] 
					\begin{tabular}{lll}

\begin{tabular}{|l|}\hline
\sf   Ahmes's selection \\ [0.01in] \hline
 $2/55=1/30+1/330\,_{\textcolor{red}{ {6}}}$ \scriptsize $\,{Eg}\normalsize$   \\ [0.01in] \hline
\end{tabular}
					 & 
\begin{tabular}{||c||c|}\hline
\sf Mother-row &  \it Divisor   \\ [0.01in] \hline
$2/11=1/6+1/66\,_{\textcolor{red}{ 6}}$ &   $\mathbf {5} $  \\ [0.01in]\hline
\end{tabular}
					&
\begin{tabular}{|c|}
\hline
\it Number of times used \\ [0.01in] \hline
$1$   \\ [0.01in] \hline
\end{tabular}
					\end{tabular}
\end{table}

Since the last case $\sf 2/95$ has no solution into 2 terms a decomposition should be found elsewhere. Where ? May be in [3-terms]  ? 
{\sf A solution within easy reach is simple} by observing that
$95=5\times 19$. $2/19$ has been \it already {\bf \cite {Brehamet1}} \rm  calculated in [3-terms] (but not $2/5$), then it is enough to divide $2/19$ by $5$.
\[
2/19=1/12+1/76\,_{\textcolor{red}{ 4}}+1/114\,_{\textcolor{red}{ 6}}.
\]
Division by $5$ yields the final result 
\begin{tabular}{|c|}\hline
$2/95=1/60+1/380\,_{\textcolor{red}{ 4}}+1/570\,_{\textcolor{red}{6}} $\scriptsize $\;{Eg}\normalsize$    \\ \hline
\end{tabular}$\;$.\\ \vspace{0.5em}
This does not undermine the role of the [2-terms] mother-table and (what is more) agrees with 
the decision  $\boxed{\sf\gamma\rm}$ discarding  any multiplicity of the last denominator beyond of $6$.\\ 
Remark: with their experience in the [4-terms] analysis and using the same methods applied to a composite number, the scribe team could have tried to search for a solution into 4 terms. 
Unfortunately the least bad solution is $2/95=1/90+1/190\,_{\textcolor{red}{ 2}}+1/285\,_{\textcolor{red}{ 3}}+1/855\,_{\textcolor{red}{ 9}}$. Of course rejected! 

									\section {Conclusion}
As we saw, Table  \ref{Composites} ($2/D$ with $D$ composite) was easily re-constructed using a global approach. Furthermore
we deliberately used not too much  arithmetical analysis because maybe that was not needed. 
We never needed assumption as to favor the largest denominator ${\mathcal D1}$. Moreover, this occurrence appears only a few times (2) in our lists. 
Used by some authors, it is not a selection criterion in our approach.
Looking  back at the successive steps of selection, it could appear as incoherent for us (modern mathematicians) because 
we would have had the reflex to analyze all case by case and in ascending order as  2/9, 2/15, 2/21 $\cdots$. 
Selections were made in  the following order:\\ 
\sf \{2/9, 2/21, 2/27, 2/33, 2/35, 2/39, 2/51, 2/57, 2/69, 2/81, 2/87, 2/93, 2/25, 2/65, 2/85, 2/49\},\\ 
\{2/15, 2/45, 2/63, 2/75, 2/99, 2/91\}, \{2/77\}, \{2/55\} and \{2/95\}.
\\ \rm
Nevertheless the choices made by the ``\it builders\rm"  of  Table  \ref{Composites} are fully coherent. \\
Actually, for composite numbers,  this table  exhibits no singularity. Singularities exist only for us.\\
{\sf All that is in favor for the fact that Egyptians have calculated all  possible cases and analyzed these ones from preliminary tables for making their choices.}


							\section*{ {\sf APPENDIX  A} : A key-equation}

\setcounter{section}{1}
\setcounter{equation}{0}

\renewcommand{\theequation }{\Alph{section}.\arabic{equation}}

\hspace{1.5em}Anyway, all the solutions for the [2-terms] cases could well have been discovered  by a method far from  our \mathversion{bold}$N^{[2]}\, \top \mbox{\tt heorem.}$\mathversion{normal}
The absence of an equational formalism in Ancient Egypt does not forbid an heuristic discovery of a key equation like equality (\ref{eq:KeyEquation}).\\ 
$D$ is decomposed as a product of two odd numbers $l$ and $m$:
\begin{equation}
D=l \times m ,\mbox{\hspace{1em}with $l \leq m$ }.
\end{equation}
Notice that this decomposition is not necessarily unique. $l,m$ may be prime numbers themselves or  between them, or not. Case $l=1$ and $m=D$  also is possible.\\

Clearly the equation used for 2-terms and $D$ prime cannot longer be exploited under a primary form analog to those used in our previous paper {\bf \cite{Brehamet1}}
for [3-terms] or [4-terms] series, namely
\begin{equation}
\mathbf{1}= \frac{D}{2D_1}+ \frac{d_2}{2D_1}.
\end{equation}
It should be transformed into a form approximately similar, by naturally imagining to\\
{\sf redistribute}   both factors $l,m$, each for one fraction (in the sense used by Bruins {\bf \cite{Bruins}}) like as:
\begin{equation}
\mathbf{1}= \frac{m}{2D_1}+ \frac{l}{2D_2},
\label{eq:UnityDec2terms}
\end{equation}
whence
\begin{equation}
\frac{2}{D}= \frac{1}{l D_1}+ \frac{1}{m D_2}.
\label{eq:InitialProblemBis}
\end{equation}
A mixed form (simultaneously additive and multiplicative) shows what is the problem to solve:
\begin{equation}
2D_1D_2= l D_1+m D_2.
\end{equation}
After a lot of trials, it depends on the intuition ability of scribes (or ourselves) for realizing that 
solutions $D_1$ could be 'centered' around $m/2$ as well as $D_2$ around $l/2$, as suggested by Eq.(\ref{eq:UnityDec2terms}).
Thus,  a first idea is to set 
\begin{equation}
2D_1=m+ r_1 \mbox{\hspace{0.5em}and \hspace{0.5em}}2D_2= l + s_2,
\end{equation}
where $r_1,  s_2$ obviously are unknown odd numbers to be found by avoiding $r_1=m$ and $s_2=l$ which would lead to $\mathcal D_1=\mathcal D_2$. \\
From the two above equations, it can derived a \sf key-equation\it,  only multiplicative\rm, and  more or less unexpected, namely 
\begin{equation}
\boxed{r_1s_2= D} \;(=lm ).
\label{eq:KeyEquation}
\end{equation}
Once a doublet  \{$r_1,s_2$\} chosen from this equation amazingly simple, \sf all the solutions of the initial problem \rm 
[see Eqs.(\ref{eq:InitialProblem}),(\ref{eq:InitialProblemBis})] \sf are available, \rm  since \sl{(in modern algebraic notations)}\rm
\begin{equation}
\mathcal D_1 =l\frac{(m+r_1)}{2}\mbox{\hspace{0.5em}and \hspace{0.5em}}\mathcal D_2 =m\frac{(l+s_2)}{2}.
\end{equation}
We cannot argue that the solutions were found by this way, of course.

							\section*{ {\sf  APPENDIX  B} : Still another approach}

\setcounter{section}{2}
\setcounter{equation}{0}

\renewcommand{\theequation }{\Alph{section}.\arabic{equation}}
From basic equation
\begin{equation}
\frac{2}{D}= \frac{1}{\mathcal D_1}+ \frac{1}{\mathcal D_2},
\label{eq:basicBis}
\end{equation}
we could search for a decomposition of $2/D$ as a sum of two other terms, each smaller than $2/D$, such that the sum to be $2/D$.
\begin{equation}
\frac{2}{D}= \frac{2}{D_1}+ \frac{2}{D_2},
\end{equation}
\it ie \rm with a ${ D_1}$ not too close to $D$ by an amount of $k_1$ and a ${D_2}$ not too close to $D$ by an amount \\of $k_2 \;(>k_1)$.
Actually this yields a decomposition of $\mathbf{1}$:
\begin{equation}
\mathbf{1}=\frac{D}{(D+k_1)}+ \frac{D}{(D+k_2)}.
\label{eq:basicTransform}
\end{equation}
With these definitions, Eq.(\ref{eq:basicTransform}) leads to 
\begin{equation}
\boxed{k_1k_2=D^2}.
\label{Babylonian}
\end{equation}
This gives the final solutions \sl{(in modern algebraic notations)}\rm
\begin{equation}
\mathcal D_1= \frac{D+k_1}{2}\mbox{\hspace{0.5em} and  \hspace{0.5em}}\mathcal D_2= \frac{D+k_2}{2}.
\end{equation}
Unknowns $k_1$, $k_2$ are not difficult to be found, with $k_1<k_2$ and $k_1$, $k_2$ odd numbers.\\
Eq.(\ref {Babylonian}) would tend to suggest  rather a Babylonian similarity, due to their techniques of using tables of squares for calculate the product of two numbers by means of mathematical identities! Besides, Ahmes says in the colophon at the beginning of papyrus that he was the copyist of an older papyrus, of Babylonian source {\bf \cite{Peet}}.  Unexpectedly this would drive us to Pythagorician triplets by the study of
\begin{equation}
(k_1+k_2)^2= (k_1-k_2)^2 +(2D)^2.
\end{equation}

\it Normally the results found in Appendix A or B lead to the same arithmetical conclusions as those obtained
by  our \mathversion{bold}$N^{[2]}\, \top \mbox{\tt heorem.}$\mathversion{normal} However it is beyond the scope of this paper to develop here these assertions.\rm

\end{flushleft}


\end{document}